\baselineskip=14pt
\font\Bf=cmbx12
\def\epsi{\varepsilon}
\def\x{x^*}
\def\y{y^*}
\def\z{z^*}
\def\u{u^*}

\def\X{X^*}

\def\<{\langle}
\def\>{\rangle}
\def\tto{\,\scriptstyle$\lower1pt\vbox{\hbox{$\to$}\kern-10pt\hbox{$\to$}}$\,}

\def\co{{\rm co}\,}
\def\dom{{\rm dom}\,}
\def\int{{\rm int}\,}
\def\implies{\ ${}\Rightarrow{}$\ }
\def\ds{\displaystyle}
\def\gata{\hfill \vrule height 6pt width 6pt depth 0pt}
\font\smc=cmcsc10
\font\srm=cmr8
\def\meaddress{\srm Department of Mathematics, University of Southern California/retired; email: verona@usc.edu}
\def\aaddress{\srm Department of Mathematics, CalState Los Angeles/retired; email: averona@calstatela.edu}
\centerline{\Bf On the Regularity of Maximal Monotone Operators and Related Results}

\bigskip
\centerline{\smc Maria Elena Verona\footnote *{\meaddress} and Andrei Verona\footnote {**}{\aaddress}}

\bigskip
\centerline{December 2012}

\bigskip
\centerline{\bf Abstract}

\bigskip
{\narrower
In the first part of the note we prove that a sufficient condition (due to Simons) for the convexity of the closure of the domain/range of a monotone operator is also necessary when the operator has bounded domain and is maximal. Simons' condition is closely related to the notion of regular maximal monotone operator. In the second part of the note we give several characterizations for the regularity of a maximal monotone operator, show that a maximal monotone operator of type (FPV) is regular and improve a previous sum theorem type result.

                                                                                                                                                                                                                                                                                                                                                                                                                                                                                                                                                                                                                                                                          } \bigskip\noindent

\bigskip\noindent
{\bf 2010 Mathematics Subject Classification:}

\noindent
Primary 47H05; Secondary 47A65.


\bigskip\noindent
{\Bf Introduction}

\bigskip\noindent
Throughout this note $X$ will denote a non zero Banach space with topological dual $\X$. For
$(x,\x)\in X\times \X$, $\<x,\x\>$ (or $\<\x,x\>$) will denote the natural evaluation map.
 If $T:X\tto \X$ is a (multivalued) operator then $G(T)=\{(x,\x)\in X\times \X: \x\in T(x)\}$ is called its {\it graph}, $D_T=\{x\in X:T(x)\ne \emptyset\}$ is called its {\it domain} and $R_T=\{\x\in\X; {\rm \ there\ exists\ }x\in X {\rm \ such\ that\ } \x\in T(x)\}$ is called its {\it range}. The operator $T$ is called {\it monotone} if $\<\x-\y,x-y\>\ge 0$ for all $(x,\x),(y,\y)\in G(T)$; $T$ is called {\it maximal monotone} if there exists no monotone operator $S$ such that $G(S)$ strictly contains $G(T)$. A pair $(x,\x)\in X\times\X$ is called {\it monotonically related to} the monotone operator $T$ if $\<\x-\z,x-z\>\ge 0$ for any $(z,\z)\in G(T)$.

Finally, throughout this note we shall refer to various well known classes of maximal monotone operators without giving their definitions (when not really necessary; this would make this note too lengthy). The reader may find the definitions in [8].

\bigskip\noindent
{\Bf 1. On the convexity of the closure of the domain of a maximal monotone operator}

\bigskip\noindent
Recall that the {\it Fitzpatrick function} $\varphi_T:X\times \X\to (-\infty,\infty]$  of a monotone operator $T$ is defined by
$$\varphi_T(x,\x)=\sup\left\{\<\x-\z,z-x\>;(z,\z)\in G(T)\right\}+\<\x,x\>.$$
It is a proper, lower semicontinuos convex function. The following statements follow immediately from the definition:

\medskip
\centerline{$(x,\x)$ is monotonically related to $T$ if and only if $\varphi_T(x,\x)\le \<x,\x\>$;}

\medskip
\centerline{If $x\in D_T$ then $\varphi_T(x,\x)\ge \<\x,x\>$ for any $\x\in\X$;}

\medskip\noindent
and

\medskip
\centerline{(1.1)\hfill $D_T\subset \co D_T\subset \pi_1\dom \varphi_T$.\hfill}

\medskip\noindent
where {\it co} represents the convex hull and $\pi_1:X\times\X\to X$ is the projection.

One of the outstanding open problems in convex analysis is whether the closure of the domain of a maximal monotone operator in an arbitrary Banach space is convex. Rockafellar [7] showed that this is true when the domain of the operator has non-empty interior. A summary of what is known about this problem can be found in [2]. In [8, Theorem 27.5] Simons proved the following result:

\bigskip\noindent
{\bf Fact 1.1.} Let $X$ be a Banach space and $T:X\tto \X $ be a monotone operator satisfying

\medskip\noindent
(1.2)\qquad $x\notin\overline D_T$ implies that $\displaystyle\sup \left\{{\<z^*, x - z\> \over ||x - z||};\, z^* \in T(z) \right \} = \infty$.

\medskip\noindent
Then:

\medskip
\item{(a)}\qquad $\pi_1 \dom \varphi_T \subset  \overline{D_T}$.

\medskip
\item{(b)}\qquad $\overline{D_T} = \overline{\co D_T} = \overline {\pi_1 \dom \varphi_T}$.

\medskip
\item{(c)}\qquad $\overline D_T$ is convex.

\bigskip\noindent
{\bf Remark.} In fact, part (a) is the important one; part (b) follows immediately from it and (1.1). Obviously (b) implies (c).

\bigskip
Recently Borwein and Yao [2, Theorem 3.6] proved the following result:

\bigskip\noindent
{\bf Fact 1.2.} Let $T:X\tto\X$ be a maximal monotone operator. Then $\overline{\co D_T} = \overline {\pi_1 \dom \varphi_T}$.

\bigskip
We shall use these results to prove a converse to Simmons' result in the case when $D_T$ is bounded and $T$ is maximal.

\bigskip\noindent
{\bf Theorem 1.3.} Let $X$ be a Banach space and  $T:X\tto \X$ be a monotone operator with bounded domain. Consider the following conditions:

\medskip
\item{(i)}\qquad $ x\notin\overline D_T$ implies that $\displaystyle\sup \left\{{\<z^*, x - z\> \over ||x - z||};\, z^* \in T(z) \right \} = \infty$.

\medskip
\item{(ii)}\qquad  $\pi_1\dom\varphi_T \subset \overline{D_T}$.

\medskip
\item{(iii)}\qquad  $\overline{D_T} = \overline{\co D_T} = \overline {\pi_1 \dom \varphi_T}$.

\medskip
\item{(iv)}\qquad  $\overline D_T$ is convex.

\medskip\noindent
Then the first three conditions are equivalent and imply the fourth one. If $T$ is maximal monotone all four conditions are equivalent.

\medskip\noindent
{\it Proof.} In view of Fact 1.1, we only need to prove that (iii) implies (i) (or (iv) implies (i) when $T$ is maximal monotone). To this end let $x\notin\overline {D_T}$ and assume  that
$$\sup \left\{{\<z^*, x - z\> \over ||x - z||};\, z^* \in T(z) \right \} < \infty.$$
It follows that there exists $M > 0$ such that $\<z^*, x - z\> \le M\|x - z\|$
for any $(z, z^*) \in G(T)$. Since $D_T$ is bounded,
$$\varphi_T ( x, 0^*) = \sup \{\< \z, x - z \>;(z, z^*) \in G(T)\}\le M \sup_{z\in D_T} \|x - z\| < \infty,$$
hence $x\in\pi_1\dom\varphi_T$. If (iii) is true, $x\in\overline{D_T}$ and this contradiction proves that (iii) implies (i). If $T$ is maximal monotone, from Fact 1.2 it follows again that $x\in \overline{\co D_T}$. Since $\overline D_T$ is convex,
$\overline{\co D_T}=\overline D_T$ and therefore $x\in \overline D_T$. This contradiction shows that (iv) implies (i).\gata

\bigskip
In [8, Theorem 27.6] Simons also proved a result dual to Theorem 1.1.

\bigskip\noindent
{\bf Fact 1.4.} Let $X$ be a Banach space and $T:X\tto \X $ be a monotone operator satisfying

\medskip\noindent
(1.3) \qquad $ \x\notin\overline R_T$ implies that $\displaystyle\sup \left\{{\<\x-\z, z\> \over \|\x - \z\|};\, z^* \in T(z)\right\} = \infty$.

\medskip\noindent
Then:

\medskip
\item{(a)}\qquad ${\pi_2 \dom \varphi_T} \subset  \overline{R_T}$.

\medskip
\item{(b)}\qquad $\overline{R_T} = \overline{\co R_T} = \overline {\pi_2 \dom \varphi_T}$.

\medskip
\item{(c)} \qquad  $\overline R_T$ is convex.

\medskip
In [3, Theorem 4.23] Borwein and Yao proved a result dual to Theorem 1.2. Here it is:

\bigskip\noindent
{\bf Fact 1.5.} Let $T:X\tto\X$ be a maximal monotone operator. Then $\overline{\co R_T}^{w^*} = \overline {\pi_2 \dom \varphi_T}^{w^*}$, where ${w^*}$ indicates the weak*-closure.

\bigskip
Finally, Theorem 1.3 has a dual counterpart which we now state.

\bigskip\noindent
{\bf Theorem 1.6.} Let $X$ be a Banach space and  $T:X\tto \X$ be a monotone  operator with bounded
range. Consider the following conditions:

\bigskip\noindent
\item{(i)}\qquad $ \x\notin\overline{R_T}$ implies that $\displaystyle\sup \left\{{\<\x-\z,z\> \over \|\x - \z\|};\, \z\in T(z)\right \} = \infty$.

\medskip\noindent
\item{(ii)}\qquad $\pi_2\dom\varphi_T \subset \overline{R_T}$.

\medskip\noindent
\item{(iii)} \qquad $\overline{R_T} = \overline{\co R_T} = \overline {\pi_2 \dom \varphi_T}$.

\medskip\noindent
\item{(iv)} \qquad  $\overline{R_T}$ is convex.

\medskip\noindent
\item{(i$'$)}\qquad $ \x\notin\overline{R_T}^{w^*}$ implies that $\displaystyle\sup \left\{{\<\x-\z,z\> \over \|\x - \z\|};\, \z\in T(z)\right \} = \infty$.

\medskip\noindent
\item{(ii$'$)} \qquad $\pi_2\dom\varphi_T \subset \overline{R_T}^{w^*}$.

\medskip\noindent
\item{(iii$'$)} \qquad $\overline{R_T}^{w^*} = \overline{\co R_T}^{w^*} = \overline {\pi_2 \dom \varphi_T}^{w^*}$.

\medskip\noindent
\item{(iv$'$)} \qquad  $\overline{R_T}^{w^*}$ is convex.

\medskip\noindent
Then (i), (ii) and (iii) are equivalent; (iii)\implies (iv)\implies (iv$'$); when $T$ is maximal monotone, (iv$'$)\implies (i$'$). In addition (i)\implies (i$'$), (ii)\implies (ii$'$)\implies (iii$'$)\implies (i$'$) and (iii$'$)\implies (iv$'$).

\medskip\noindent
{\it Proof.} From Fact 1.4 it follows that (i)\implies (ii)\implies (iii)\implies (iv). Exactly as in the proof of Theorem 1.3, we can prove that (iii)\implies (i) and (iii$'$)\implies (i$'$). Clearly (i)\implies (i$'$), (ii)\implies (ii$'$)\implies (iii$'$)\implies (iv$'$). When $T$ is maximal monotone, exactly as in the proof of (iv)\implies (i) in Theorem 1.3 (but using Fact 1.5 instead of Fact 1.2), we can show that (iv$'$)\implies (i$'$). It remains to show that (iv)\implies (iv$'$). So assume that $\overline{R_T}$ is convex. Then
$$\co R_T\subset\overline{R_T}\subset \overline{R_T}^{w^*}$$
and it follows that
$$\overline{R_T}^{w^*}\subset\overline{\co R_T}^{w^*}\subset\overline{(\overline{R_T})}^{w^*}\subset \overline{R_T}^{w^*}.$$
Therefore $\overline{R_T}^{w^*}=\overline{\co R_T}^{w^*}$, showing that $\overline{R_T}^{w^*}$ is convex.\gata

\bigskip\noindent
{\Bf 2. Regular maximal monotone operators}

\bigskip\noindent
In [10] we studied monotone operators with properties similar to (1.2) and introduced the following number
$$L(x,\x,T)=\sup\left\{{\<\y-\x,x-y\>\over \|x-y\|}:\,y\ne x, (y,\y)\in G(T)\right\} $$
With this notation, condition (1.2) becomes
$$ L(x,0^*,T) < \infty \Rightarrow x \in \overline{D_T}.$$
We turn now our attention to monotone operators with a more restrictive property than (1.2), namely monotone operators with the property:

\medskip
\centerline{$ L(x,0^*,T) < \infty \Rightarrow x \in D_T$ ({\it regularity property}).}

\medskip\noindent
It follows immediately that if $T:X\tto \X$ is a monotone operator  satisfying the regularity property then
$\pi_1 \dom \varphi_T \subset  \overline{D_T}$, hence $ \overline{D_T} = \overline{\co D_T} = \overline {\pi_1 \dom \varphi_T}$ and $\overline D_T$ is convex. As a matter of fact, we also proved this fact in [10] by using different methods.

\medskip\noindent
Recall that (see [10]) $$L(x,\x,T)\le \inf\left\{\|\y-\x\|: \y\in T(x)\right\}=d(\x,T(x))$$

\medskip\noindent
where $d$ represents the distance from a point to a set.

\bigskip\noindent
{\bf Definition 2.1.} (see [10]) A maximal monotone operator $T$ is called {\it regular} if

\medskip
\centerline {$L(x,\x,T)=d(\x,T(x))$ for any $(x,\x)\in X\times \X$.}

\bigskip\noindent
{\bf Remark 2.2.}   If for $(x, x_0^*)\in X\times \X$, $L(x, x_0^*, T) < \infty$ then $L(x,x^*, T) < \infty$ for all $x^* \in X^*$.  Indeed, for any $(\y,y)\in G(T)$ we have
$${\<\x-\y,y-x\>\over \|x-y\|}\le{\<\x_0-\y,y-x\>\over \|x-y\|}+{\<\x-\x_0,y-x\>\over \|x-y\|}\le {\<\x_0-\y,y-x\>\over \|x-y\|}+\|\x-\x_0\|.$$

The next result contains some known characterizations of regularity and also some new ones. Before we state it, recall that for any $\lambda >0$ and $x\in X$, $g_{\lambda,x}:X\to R$ is defined by $g_{\lambda,x}(y)=\lambda\|x-y\|$.

\bigskip\noindent
{\bf Theorem 2.3.} Let $T:X\tto \X$ be a maximal monotone operator on a Banach space $X$.
The following are equivalent:

\medskip
\item{(a)} \quad $L(x,\x, T) = d(x^*,T(x))$ for all $(x,\x) \in X \times X^*$ (i.e. $T$ is regular).

\medskip
\item{(b)} \quad If $S:X\tto \X$ is a bounded maximal monotone operator then $T+S$ is maximal monotone.

\medskip
\item{(c)} \quad $T+\partial g_{\lambda, x}$ is maximal monotone for all $\lambda > 0$ and $x\in X$.

\medskip
\item{(d)} \quad If $L(x,0^*,T)<\infty$ then $x\in D_T$ (i.e. $T$ has the {\it regularity property)}.

\medskip
\item{(e)} \quad $L(x,0^*,T)= d(0^*,T(x))$ for all $x\in X$.

\medskip\noindent
{\it Proof.} (a) $\Rightarrow$ (b) was proved in [12, Theorem 1.1], clearly (b) $\Rightarrow$ (c) while (c) $\Rightarrow$ (a) was proved in [10]. It is also obvious that (a) $\Rightarrow$ (d). Next assume that (d) is true. We shall show that (c) is true. In view of [13, Theorem 3.4 and Corollary 5.6] or [8, Theorem 24.1] it is enough to check that $\varphi_{T+\partial g_{\lambda,x}}(x,x^*)\ge\<x,\x\>$ on $X\times X^*$. To this end let $(x, x^*)$ be such that $\varphi_{T+\partial g_{\lambda,x}}(x,x^*)\le\<x,\x\>$, i.e. $(x,\x)$ is monotonically related to $T + \partial g_{\lambda,x}$. An easy computation shows that $L(x,\x, T)\le\lambda$, hence $L(x,0^*,T)<\infty$. Thus $x\in D_T=D_{T+\partial g_{\lambda,x}}$. Then $\varphi_{T+\partial g_{\lambda,x}}(x,\x)\ge \<x,\x\>$ and, as pointed above, $T+\partial g_{\lambda, x}$ is maximal monotone. Obviously (a) $\Rightarrow$ (e) and (e) $\Rightarrow$ (d). The theorem is proved.\gata

\bigskip\noindent
{\bf Remark 2.4.} Any monotone operator that has the regularity property and with $T(x)$ closed for $x \in D_T$ is a maximal monotone operator  (see [10], Remark 3).

\bigskip\noindent
{\bf Proposition 2.5.} A maximal monotone operator $T:X\tto \X$ is regular if and only if $T+S$ is a regular maximal monotone operator for any bounded maximal monotone operator $S:X\tto \X$.

\medskip\noindent
{\it Proof.} Assume first that $T$ is regular. Let $U:X\tto \X$ be a bounded maximal monotone operator. Being maximal and bounded, both $S$ and $U$ must have full domain (see for example [8, Theorem 25.1]) and in view of Heisler's Theorem (see [6, Remark on page 17]) $S+U$ is maximal monotone, obviously bounded. Theorem 2.3 implies first that $T+S+U$ is maximal monotone and second that $T+S$ is regular. Obviously, if $T + S$ is a regular maximal monotone operator for any bounded maximal monotone operator $S$, by taking $S = 0$, we obtain that $T$ is a regular maximal monotone operator.\gata

\medskip
It is well known that if $T:X\tto \X $ is a maximal monotone operator on a reflexive Banach space $X$, then $\overline{D_T} = \overline{\co D_T} = \overline {\pi_1 \dom \varphi_T}$. This result remains true for maximal monotone operators of type (FPV) in arbitrary Banach spaces [8, Theorem 44.2]. From the discussions on monotone operators with the regularity property it follows that in arbitrary Banach spaces the result remains true
for regular maximal monotone operators. As a matter of fact, we shall prove that a maximal monotone operator of type (FPV) is regular.

\medskip\noindent
Here is a list of maximal monotone operators that are regular:

\medskip
\item{-} subdifferentials of proper, convex, lower semi-continuous functions (they satisfy (c) in Theorem 2.3);

\item{-} maximal monotone operators on reflexive Banach spaces [10];

\item{-} strongly-representable maximal monotone operators (or equivalently maximal monotone operators of type (NI), see [1]) are regular [15, Remark 7]; in particular any maximal monotone operator of type (D) is regular.

\item{-} linear (possibly discontinuous) maximal monotone operators [11, Theorem 3.2];

\item{-} strongly maximal monotone operators [10, Proposition 1];

\item{-} maximal monotone operators with closed and convex domains (follows from [13, Theorem 5.10 ($\beta$)]);

\item{-} maximal monotone operators such that their domain has non-empty interior [12, Corollay 1.3].

\medskip
Next, we shall show that the class of monotone operators of type (FPV) is included in the class of regular maximal monotone operators.

\bigskip\noindent
{\bf Theorem 2.6.} Let $T:X\tto \X $ be a maximal monotone operator of type (FPV), then $T$ is a regular maximal monotone operator.

\medskip\noindent
{\it Proof.} According to [16, Corollary 3.6] the sum of a maximal monotone operator of type (FPV) with a maximal monotone operator with full domain is maximal. It follows that $T + \partial g _{\lambda, x}$ is maximal monotone for all $\lambda \ge 0$ and $x\in X$; hence $T$ is a regular maximal monotone operator by Theorem 2.3.\gata

\bigskip\noindent
{\bf Definitions 2.7.}
For a monotone operator $T:X\tto X^*$ we recall some notation and definitions from in [11] and introduce some new ones:

\medskip
\centerline{${\cal C}_0(T) = \bigl\{C \subset X;\, C {\rm\ is \ closed, \ convex, \ and \ } D_T \cap \int (C) \ne \emptyset \bigr\}$.}

\medskip
\centerline {${\cal C}_1(T) = \bigl\{ C \subset X;\, C {\rm\ is \ closed, \ convex, \ and \ } \bigcup_{\lambda > 0} \lambda (\co (D_T) - C)$ is a closed subspace of $X \bigr\}.$}

\medskip
$T$ is called ${\cal C}_i$-{\it maximal} if $T+\partial I_C$ is a maximal monotone for any $C\in {\cal C}_i(T)$, $i\in\{0,1\}$.

\medskip
$T$ is called ${\cal C}_i$-{\it regular} if $T+\partial I_C$ is a regular maximal monotone for any $C\in {\cal C}_i(T)$, $i\in\{0,1\}$.

\bigskip\noindent
{\bf Remarks.} (1) Let $T$ be ${\cal C}_i$-maximal (respectively regular), $i\in \{0,1\}$. Since $X\in {\cal C}_i(T)$ then $T$ is maximal (respectively regular).

\medskip
(2) If $T$ is ${\cal C}_1$-maximal (respectively regular), then $T$ is ${\cal C}_0$-maximal (respectively regular).

\medskip
 We show next that certain classes of maximal monotone operators are ${\cal C}_0$ or ${\cal C}_1$-regular.

\bigskip\noindent
{\bf Proposition 2.8.} Let $T:X\tto\X$  be a maximal monotone operator.

\medskip
\item{(1)} If $D_T$ is closed convex (in particular $D_T = X$) and $C \in {\cal C}_1 (T)$, then $T + \partial I_C$ is a regular maximal monotone operator. As a consequence, $T$ is a ${\cal C}_1$-regular maximal monotone operator.

\medskip
\item{(2)} If $T$ is (FPV) with $D(T)$ convex and $C \in {\cal C}_0(T) $, then $T + \partial I_C$ is a  regular maximal monotone operator. As a consequence, $T$ is a ${\cal C}_0$-regular maximal monotone operator.

\medskip\noindent
{\it Proof.} (1) If $D_T$ is closed convex it follows from [13, Theorem 5.10 ($\beta$)] that $T +\partial I_C$ is maximal monotone, hence by the same theorem
$T + \partial I_C + \partial g _{\lambda, x}$ is a maximal monotone operator for any $\lambda >0$ and any $x\in X$. From Theorem 2.3 it follows that $T + \partial I_C$ is a regular maximal monotone
operator, hence $T$ is a ${\cal C}_1$-regular maximal monotone operator. The particular case when $D_T = X$ has an earlier proof [11, Corollary 2.11(i)].

(2) Under these assumptions, it follows from [16, Corollary 2.10] that $T + \partial I_C$, is of type (FPV). Theorem 2.6 implies that $T+\partial I_C$ is a regular maximal monotone operator.\gata

\bigskip\noindent
{\bf Lemma 2.9.} Assume that $D,C,K\subset X$, ${\overline D}, C, K$ are closed convex, $D\cap\int C\ne\emptyset$ and $D\cap C\cap\int K\ne\emptyset$. Then $D\cap\int C\cap\int K\ne\emptyset$.

\medskip\noindent
{\it Proof}. Let $x\in D\cap\int C$ and $y\in D\cap C\cap\int K$. Then $[x,y]\subset \overline D$, $[x,y)\subset \int C$ and there exists $z\in (x,y)\cap K$ such that $[z,y]\subset \int K$. We can choose $z$ sufficiently close to $y$ such that $z\in\overline D\cap\int C\cap\int K$. This implies that $D\cap\int C\cap\int K\ne\emptyset$.\gata

\bigskip\noindent
{\bf Lemma 2.10.} Assume that $C\in{\cal C}_i(T)$, $i\in\{0,1\}$, and $K\in{\cal C}_0(T+\partial I_C)$; if $i=0$ assume also $\overline D_T$ is convex. Then $C\cap K\in {\cal C}_i(T)$.

\medskip\noindent
{\it Proof}. If $i=0$, the statement follows from the previous lemma. Let $i=1$. We can and do assume that $0\in D_T\cap C\cap \int K$. Let $x\in Y=\bigcup_{\lambda>0}(\co D_T-C)$. Since $C\in{\cal C}_1(T)$, $Y$ is closed and there exists $s>0$ such that $x=s(y-z)$ with $y\in\co D_T$ and $z\in C$. We can choose $t>0$, sufficiently small such that $tz\in K$; clearly $ty\in \co D_T$ and $tz\in C$. It follows that $x={s\over t}(ty-tz)$ and therefore $Y\subset \bigcup_{\lambda>0}(\co D_T-C\cap K)$.
Finally, since $\bigcup_{\lambda>0}(\co D_T-C\cap K)\subset \bigcup_{\lambda>0}(\co D_T-C)=Y$ the lemma is proved.\gata

\bigskip\noindent
{\bf Theorem 2.11.} Let $T:X\tto\X$ be ${\cal C}_i$-maximal ($i=0,1$) and let $C\in{\cal C}_i(T)$. Then $T+\partial I_C$ is of type (FPV). In particular $T+\partial I_C$ is regular and thus $T$ is $C_i$-regular.

\medskip\noindent
{\it Proof}. Let $i = 0$ and $T$ be ${\cal C}_0$-maximal. It follows from Lemma 2.10 that $T + \partial I_C$ is ${\cal C}_0$-maximal, hence $T + \partial I_C$ is of type (FPV), implying that $T + \partial I_C$ is regular. Thus $T$ is ${\cal C}_0$-regular. In fact, more generally, $T+\partial I_C$ is ${\cal C}_0$-regular.

\medskip
Let $i=1$ and $T$ be ${\cal C}_1$-maximal. Take $C\in {\cal C}_1(T)$ and $K\in {\cal C}_0(T + \partial I_C)$. Since $T$ is of type (FPV), $T$ is regular, hence $\overline {D_T}$ is convex. Hence $C \cap K \in {\cal C}_1(T)$ (see Lemma 2.10). It follows that $T + \partial I_{C \cap K} = T + \partial I_C + \partial I_K$ is maximal. This implies that $T + \partial I_C$ is of type (FPV), therefore $T + \partial I_C$ is regular. Thus $T$ is ${\cal C}_1$-regular.\gata

\bigskip\noindent
{\bf Notation}. If $x\in X$ and $A\subset X$ set $K_{x,A} = \big \{(1-t)x + ta ;\, 0\le t\le 1,a\in A\big \}$.

\bigskip\noindent
{\bf Definition 2.12.} A monotone operator $T$ is called {\it directionally inf bounded at} $x\in D_T$ if for any $y\in X$ there exists $\epsilon, M>0$ such that for any $z\in K_{y,B(x,\epsilon)}\cap D_T$, there exists $z^*\in T(z)$ with $\<z^*,z -y\>\le M\|z-y\|$.

\medskip\noindent
{\bf Lemma 2.13.}. If $T:X\tto X^*$ is maximal monotone and $\int D_T\ne\emptyset$ then $T$ is directionally inf bounded at any $x\in \int(D_T)$

\medskip\noindent
{\it Proof.} Let $x\in\int(D_T)$. It is well known [7] that $T$ is bounded at such a point and thus there exists $\epsi>0$ and $M>0$ such that $B(x,\epsilon)\subset D_T$ and for any $u\in B(x,\epsi)$ and $\u\in T(u)$ we have $\|\u\|\le M$.
Let $y\in X$ and $z\in K_{y,B(x,\epsilon)}\cap D_T$, $z\ne y$. Then there exists $u\in B(x,\epsi)$ and $t\in(0,1]$ such that $z=(1-t)y+tu$. Choose $\z\in T(z)$. If $z=u$ then obviously $\<\z,z-y\>\le M\|z-y\|$.
If $z\ne u$ then $t\ne 1$. Choose $\u\in T(u)$. An easy computation shows that $z-y={t\over 1-t}(u-z)$ and therefore

\medskip
\centerline{$\ds\hfill{\<\z,z-y\>\over\|z-y\|}={\<\z,u-z\>\over\|u-z\|}={\<\z-\u,u-z\>\over\|u-z\|}+{\< u^*,u-z\>\over\|u-z\|}\le M.\gata$}

\medskip
The following result is an expansion of our sum theorem in [11], with a simpler proof.

\bigskip\noindent
{\bf Theorem 2.14.} Let $T,\ S: X\tto X^*$ be maximal monotone operators with $D_T\cap \overline {D_S}\subset D_S$.

\medskip
\item{(1)}  If $T$ is $C_1$-maximal, ${\overline D_S}\in {\cal C}_1(T)$ and there exists $x_0\in D_T\cap D_S$ such that $S$ is directionally inf bounded at $x_0$, then $T + S$ is a maximal monotone operator.

\medskip
\item{(2)} If $T$ is $C_0$-maximal and $D_T\cap \int D_S\ne\emptyset$ then $T + S$ is a maximal monotone operator.

\medskip\noindent
{\it Proof.} (1) Assume that $(x,\x)\in X\times\X$ is monotonically related with $G(T+S)$. Set $K=K_{x,B(x_0,\epsi)}$. Then for any $z \in K\cap D_S$ there exists $z^* \in S(z)$ with $\<z^*,z-x\>\le  M\|z-x\|$. Let $z \in D_T \cap K\cap \overline{D_S}\subset D_T \cap K\cap D_S$, $t^* \in T(z) $, and $u^* \in \partial I_{K \cap\overline{D_S}}(z) = \partial (I_K + I_{\overline{D_S}}) (z)$. From the sum formula for subdifferentials, there exits $j^* \in \partial I_K (z)$, and $i^* \in \partial I_{\overline{D_S}}(z)$ such that $u^* = j^* + i^*$. We also have $\<j^*,x-z\>\le 0$ and $s^*+i^*\in S(z)$. Then
$${ \< t^* + u^* - x^*, x - z\>\over ||x - z||} = {{\<t^* + s^* + i^* - x^*, x - z\>}\over {||x - z||}} +
{\<j^*, x - z\> \over ||x - z||} - {\<s^*, x - z\> \over ||x - z||} \le 0 + M.$$
Thus $L(x,x^*,T+\partial I_{K\cap\overline {D_S}})<\infty$. Obviously ${\overline D_S}\in C_1(T)$ and by Lemma 2.10 ${\overline D_S}\cap K\in C_1(T)$. Since $T$ is ${\cal C}_1$-maximal, Theorem 2.11 implies that $T$ is ${\cal C}_1$-regular and therefore $x\in D_T\cap K\cap\overline {D_S} \subset D_T\cap D_S=D_{T+S}$, implying that $\varphi(x,\x)\ge\<x,\x\>$. This shows that $T + S$ is maximal monotone (see [13, Theorem 3.4 and Corollary 5.6] or [8, Theorem 24.1]).


\medskip
(2) This assertion can be proved in the same way as (1) or follows from [16, Corollary 3.9] since any ${\cal C}_0$-maximal monotone operator is of type (FPV).\gata

\bigskip\noindent
{\Bf References}

\medskip
\hang\noindent
[1] M. Marques Alves and B. F. Svaiter: A new old class of maximal monotone operators, J. Convex Analysis, 16 (2009), 881--890.

\medskip
\hang\noindent
[2] J. M. Borwein and L. Yao: Some results on the convexity of the closure of the domain of a maximally monotone operator, Preprint, arXiv:1205.4482v1, May 2012.

\medskip
\hang\noindent
[3] J. M. Borwein and L. Yao: Recent progress on monotone operator theory, Preprint, arXiv:1210.3401v2, October 2012.

\medskip
\hang\noindent
[4] S. P. Fitzpatrick, Representing monotone operators by convex functions, in Workshop/Miniconference on Functional Analysis and Optimization (Canberra 1988), Proc. Centre Math. Anal., Austral. Nat. Univ., vol. 20, (1988), Canberra, Australia, pp. 59--65.

\medskip
\hang\noindent
[5] S. P. Fitzpatrick and R. R. Phelps: Some properties of maximal monotone operators on nonreflexive Banach spaces, Set-Valued Anal. 3 (1995), 51--69.

\medskip
\hang\noindent
[6] R. R. Phelps: Lectures on Maximal Monotone Operators, Extracta Mathematicae, 12 (1997), 193--230.

\medskip
\hang\noindent
[7] R.T. Rockafellar: Local boundedness of nonlinear, monotone operators, Michigan Math. J. 16 (1969), 397--407.

\medskip
\hang\noindent
[8] S. Simons: From Hahn-Banach to Monotonicity, Second Edition, Lecture Notes in Mathematics 1693 (2008), Springer-Verlag.

\medskip
\hang\noindent
[9] A. Verona and M. E. Verona: Remarks on subgradients and $\epsi$-subgradients, Set-Valued Anal. 1 (1993), 261--272.

\medskip
\hang\noindent
[10] A. Verona and M. E. Verona: Regular maximal monotone operators, Set-Valued Anal. 6 (1998), 302--312.

\medskip
\hang\noindent
[11] A. Verona and M. E. Verona: Regular maximal monotone operators and the sum theorem, J. Convex Analysis, 7 (2000), 115--128.

\medskip
\hang\noindent
[12] A. Verona and M. E. Verona: Regular maximal monotone multifunctions and enlargements, J. Convex Analysis, 16 (2009) 1003--1009.

\medskip
\hang\noindent
[13] M. D. Voisei: The sum and chain rules for maximal monotone operators, Set-Valued Analysis, 16 (2008), 461--476.

\medskip
\hang\noindent
[14]: M. D.  Voisei: A sum theorem for (FPV) operators and normal cones, J. Math. Anal. Appl., vol. 371 (2010), pp. 661--664.

\medskip
\hang\noindent
[15] M. D. Voisei and C. Zalinescu: Strongly-representable monotone operators, J. Convex Analysis, 16 (2009), 1011--1033.

\medskip
\hang\noindent
[16] L. Yao: The sum of a maximal monotone operator of type (FPV) and a maximal monotone operator with full domain is maximal monotone, Nonlinear Analysis 74 (2011) 6144--6152.

\bye